\documentclass[11pt,oneside,a4paper,reqno]{amsart}
\usepackage[T1]{fontenc}
\usepackage[english]{babel}
\usepackage{amssymb,amsmath,amsthm,graphicx,amsfonts}
\usepackage[foot]{amsaddr}
\usepackage{mathrsfs}
\usepackage{lmodern}
\RequirePackage{fix-cm} 

\usepackage[dvipsnames]{xcolor}
\definecolor{ourcolor}{RGB}{66, 103, 138}

\usepackage[
   rm={oldstyle=false,proportional=true},
   sf={oldstyle=true,proportional=true},
   tt={oldstyle=true,proportional=true,variable=true},
   qt=false,
 ]{cfr-lm}
\usepackage{cite,enumerate,setspace}
\usepackage[normalem]{ulem}
\usepackage{hyperref}
 \hypersetup{
 colorlinks,
 citecolor=ourcolor,
 linkcolor=ourcolor,
 urlcolor=black}

\numberwithin{equation}{section}

\usepackage{geometry}
\let\originalleft\left
\let\originalright\right
\renewcommand{\left}{\mathopen{}\mathclose\bgroup\originalleft}
\renewcommand{\right}{\aftergroup\egroup\originalright}

 \makeatletter
 \def\@seccntformat#1{\hspace*{0mm}%
  \protect\textup{\protect\@secnumfont
    \ifnum\pdfstrcmp{subsection}{#1}=0 \bfseries\fi
    \csname the#1\endcsname
    \protect\@secnumpunct
      }%
 }

\newcommand{\assign}{:=}
\newcommand{\divides}{\mathrel{|}}
\newcommand{\mathLaplace}{\Delta}
\newcommand{\of}{:}
\newcommand{\tmabbr}[1]{#1}

\newcommand{\tmem}[1]{{\em #1\/}}
\newcommand{\tmmathbf}[1]{\ensuremath{\boldsymbol{#1}}}
\newcommand{\tmname}[1]{\textsc{#1}}
\newcommand{\tmop}[1]{\ensuremath{\operatorname{#1}}}
\newcommand{\tmrsub}[1]{\ensuremath{_{\textrm{#1}}}}
\newcommand{\tmstrong}[1]{\textbf{#1}}
\newcommand{\tmtextbf}[1]{\text{{\bfseries{#1}}}}

\newcommand{\tmtextsc}[1]{\text{{\scshape{#1}}}}

\newcommand{\tmtextup}[1]{\text{{\upshape{#1}}}}
\newenvironment{enumerateroman}{\begin{enumerate}[i.] }{\end{enumerate}}
\newtheorem{definition}{Definition}
\newtheorem{lemma}{Lemma}
{\theoremstyle{remark}\newtheorem{remark}{Remark}}
\newtheorem{theorem}{Theorem}

\newcommand{\jqed}{{\hspace*{\fill}}$\square$}
\newcommand{\gradsym}{\nabla_{\tmop{sym}}}
\newcommand{\gradskw}{\nabla_{\tmop{skw}}}
\newcommand{\divbsym}{\textrm{\tmtextbf{div}}_{\tmop{sym}}}
\newcommand{\divbskw}{\textrm{\tmtextbf{div}}_{\tmop{skw}}}
\newcommand{\divb}{\textrm{\tmtextbf{div}}}

\newcommand{\jdual}[1]{\tmmathbf{\langle} #1 \tmmathbf{\rangle}}
\newcommand{\Om}{\Omega}

\newcommand{\cpt}{\tmtextup{)}}
\newcommand{\opt}{\tmtextup{(}}
\newcommand{\Fsp}{\textbf{\,{\color[HTML]{444466}{\of}}\,}}
\newcommand{\T}{\mathsf{T}}
\newcommand{\curl}{\mathbf{curl}}
\newcommand{\lapl}{\Delta}
\newcommand{\grad}{\nabla}
\newcommand{\divv}{\text{div}}

\newcommand{\RR}{\mathbb{R}}
\newcommand{\NN}{\mathbb{N}}

\newcommand{\eqs}{=}

\begin{document}

\title[{\sc Korn and Poincar{\'e}-Korn inequalities}]{{\large \bf \uppercase{Korn and Poincar{\'e}-Korn inequalities: \\
A different perspective}}}

\author[1]{{\sc Giovanni Di Fratta}$^{(1)}$}

\author[2]{{\sc Francesco Solombrino}$^{(1)}$}

\address[1]{Dipartimento di Matematica e Applicazioni ``R.~Caccioppoli'', Università degli Studi di Napoli
    ``Federico II'', Via Cintia, 80126, Napoli, Italy.}

\begin{abstract}
  We present a concise point of view on the first and the second Korn's inequality for general exponent $p$ and for a class of domains that includes Lipschitz domains. Our argument is conceptually very simple and, for $p = 2$, uses only the classical Riesz representation theorem in Hilbert spaces. Moreover, the argument for the general exponent $1<p<\infty$ remains the same, the only change being invoking now the $q$-Riesz representation theorem (with $q$ the harmonic conjugate of $p$). We also complement the analysis with elementary derivations of Poincar{\'e}-Korn inequalities in bounded and unbounded domains, which are essential tools in showing the coercivity of variational problems of elasticity but also propedeutic to the proof of the first Korn inequality.
   \vspace{4pt}
  
  \noindent {\scriptsize {\sc Mathematics Subject Classification.} 35A23, 46N20, 46F10, 74B20.}
  
  \vspace{1pt}
  
  \noindent {\scriptsize {\sc Key words.} Korn inequality, Poincar{\'e}-Korn inequality, Riesz representation theorem, Weyl lemma.}
\end{abstract}

\begingroup
\def\uppercasenonmath#1{} 
\let\MakeUppercase\relax 
\maketitle
\endgroup

\section{Introduction and motivation}
\noindent Korn's inequality is a fundamental result in mathematical analysis and a
cornerstone of linear elasticity theory. For an elastic body occupying a
regular region $\Om \subset \RR^3$, or more generally for a smooth bounded
domain $\Om \subset \RR^N$, $N \geqslant 1$, the inequality, in its more
classical form, states the existence, for every $1 < p < \infty$, of a
positive constant $K_{\Om, p} > 0$, depending only on $\Om$ and $p$, such that
\begin{equation}
  \int_{\Om} | \tmmathbf{u} |^p + | \grad \tmmathbf{u} |^p
  \leqslant K_{\Om, p} \int_{\Om} | \tmmathbf{u} |^p + | \gradsym
  \tmmathbf{u} |^p \quad \quad \forall \tmmathbf{u} \in C^{\infty}
  (\overline{\Om}, \RR^N), \label{eq:K1W1pelementary}
\end{equation}
with $\gradsym \tmmathbf{u}$ being the symmetric part of the Jacobian matrix
$\grad \tmmathbf{u}$, $| \gradsym \tmmathbf{u}|$ and $| \grad \tmmathbf{u}|$
the corresponding Frobenius norms (see sec.~\ref{sec:notation}).

In the language of the theory of elasticity, the symmetric part of the
gradient represents the strain experienced by an elastic body when subjected
to deformation. Inequality \eqref{eq:K1W1pelementary} allows controlling the
$L^p$-norm of the gradient of the displacement field with the $L^p$-norm of
the {\tmem{linearized}} strain tensor. It is a fundamental tool for
establishing the existence, uniqueness, and stability of solutions to linear
and nonlinear elastic deformation problems. However, it has applications in
diverse areas of mathematics and physics. The family of open sets where the
inequality \eqref{eq:K1W1pelementary} holds, sometimes referred to as
{\tmem{Korn domains}}~{\cite{Acosta2017a}}, allows for a well-posed theory of
linear elasticity and, therefore, for reliable reduced models, for example,
for the study of elastic plates and
shells~{\cite{Ciarlet[2022]copyright2022b,Ciarlet[2022]copyright2022a}}.

The inequality \eqref{eq:K1W1pelementary} is usually referred to as the
{\tmem{second}} Korn
inequality~{\cite{Ciarlet[2022]copyright2022,Marsden1994,Oleinik1992}} to
distinguish it from its variant in the restricted class
$C^{\infty}_c (\Om, \RR^N)$ of vector fields with compact support
in $\Om$, usually known as the {\tmem{first}} Korn inequality, which states
the existence of a positive constant $K_{\Om, p} > 0$, such that
\begin{equation}
  \int_{\Om} | \grad \tmmathbf{u} |^p \leqslant K_{\Om, p}
  \int_{\Om} | \gradsym \tmmathbf{u} |^p \quad \quad \forall
  \tmmathbf{u} \in C^{\infty}_c (\Om, \RR^N) .
  \label{eq:K1W1pelementarycc}
\end{equation}
The first Korn inequality, for $p = 2$, immediately follows from the equality
\begin{equation}
  \int_{\Om} | \gradsym \tmmathbf{u}|^2 \eqs \frac{1}{2} \int_{\Om} | \grad
  \tmmathbf{u}|^2 + | \tmop{div} \tmmathbf{u} |^2 \quad \quad \forall
  \tmmathbf{u} \in C^{\infty}_c (\Om, \RR^N).
  \label{eq:simpleidentity}
\end{equation}
Note that, when $p=2$, the constant in \eqref{eq:K1W1pelementarycc} can be chosen independently of $\Omega$ (for example, $K_{\Omega, 2} = 2$). This turns out to be a general remark. Indeed, as long as one is not interested in optimal constants, a simple scaling argument on balls gives that the constant $K_{\Om, p}$ in equation \eqref{eq:K1W1pelementarycc} can always be taken to be independent of $\Om$.

Despite their significance, available proofs of Korn inequalities are pretty
involved. Their proof strategies typically rely on approximation and
compactness arguments combined with Calderón-Zygmund estimates for singular
integral operators. As described in subsection~\ref{sec:stateoftheart}, over
the years, there has been a sporadic yet constant scientific production of
different points of view that often led to simplifications or extensions.

The primary objective of this paper is to present a concise point of view on
the first and the second Korn's inequality, for general exponent $p$, and for
a class of domains that includes Lipschitz domains. Our proof is conceptually
very simple and, for $p = 2$, gives a proof of the second Korn's inequality
that uses only the classical Riesz representation theorem in Hilbert spaces.
Moreover, conceptually, the argument for the general exponent $p$ remains the
same, the only change being invoking now the $q$-Riesz representation theorem
(with $q$ the harmonic conjugate of $p$) instead of the classical $2$-Riesz
representation theorem.
We stress that the chief outcomes of the article are known, 
the novelty being in the alternative take on the subject we propose. 
Further details on the approach are given in section~\ref{sec:setup}.

\subsection{State of the art}\label{sec:stateoftheart}Named after the German
scientist Arthur Korn\footnote{Arthur Korn was a mathematician, physicist, and
inventor. Korn was born on May 20, 1870, in Breslau, Germany (now Wroc{\l}aw,
Poland). He studied physics and mathematics first at the University of
Freiburg/Breisgau (with Emil Warburg) and then at the University of Leipzig
(with Carl Neumann), where he graduated in 1890. Afterward, he studied in
Berlin, Paris, London, and W{\"u}rzburg, benefiting from the guidance of other
prominent figures like Henri Poincar{\'e} during his time in Paris. In 1895,
he became Priv.-Doz. in physics at the University of Munich, where he was
appointed professor in 1903. In 1914, he accepted the chair of physics at the
Berlin Institute of Technology. He emigrated to the USA with his family in
1939, where he became a professor of mathematics and physics at the Stevens
Institute of Technology in Hoboken (New Jersey) in the same year. In 1945, he
received American citizenship. He died at the age of 75 in the Jersey City
Medical Center~{\cite{NYT45}}.
To the general public, Arthur Korn is primarily recognized for his pioneering
contributions to telecommunications engineering. In 1902, he patented the
first apparatus that enabled the transmission of images through telegraph
lines, and with this system, in 1904, as part of a public demonstration, a
portrait of Kaiser Wilhelm II from Munich to Nuremberg. In 1907, The New York
Times published an article headlined ``{\tmname{Photographs by Telegraph:
Television Next?}}'' that boldly declared that the {\guillemotleft}successful
test of Prof. Korn's remarkable invention indicates the possibility of another
field for scientific discoverer{\guillemotright}~{\cite{NYT2016}}. For a
comprehensive account of Arthur Korn's life, including his groundbreaking
contributions to telephotography, we refer to~{\cite{KornStory1950}}.}, Korn
inequalities are a subject with more than a century of history. Their
extensions to more general settings continue to be an active area of research
(see, e.g.,
{\cite{Breit2012,Cagnetti2022,MR3250367,MR3096302,MR3194689,Desvillettes2002,MR3827187,Grabovsky2018,Jiang2017,MR3294348}}).
The literature on the topic is enormous, and in what follows, we only focus on
the (still vast) literature that is most pertinent to our study. For
complementary aspects of the questions, we refer the reader to the excellent
and comprehensive treatments
{\cite{Ciarlet[2022]copyright2022,Marsden1994,Oleinik1992}}, {\cite[secs.~XI.4
and XII.4]{Dautray1990}}, {\cite[sec.~ 3.3]{Duvaut1972}}, and to
{\cite[sec.~7]{Horgan1995}} for a quick survey on the applications of Korn
inequalities to continuum mechanics. As far as possible, we present the
relevant literature in historical chronological order.

The first and the second Korn's inequality, for $p = 2$, were formulated and
proved by Korn himself ({\tmabbr{cf.}}~{\cite{Korn1907,korn1909einige}}).
While the original proof in {\cite{Korn1907}} for the first case is clear and
essentially based on the identity \eqref{eq:simpleidentity}, the proof of the
second Korn inequality given by Korn in {\cite{korn1909einige}} is, in the
words of Friedrichs~{\cite{Friedrichs1947}}, {\tmem{very complicated}}, to the
level that Friedrichs himself had to admit he was unable to verify Korn's
proof for the second case. The reason for that, according to Friedrichs, other
than the length of the involved arguments, was their nested structure: the
proof given in {\cite{korn1909einige}} refers to previous results presented in
{\cite{korn1908solution}}, which still refer to earlier papers of Korn. Giving
direct proof of Korn's inequality in the second case was what motivated
Friedrichs contribution~{\cite{korn1909einige}}.

Starting from the work of Friedrichs, several authors proposed alternative
proofs and generalizations. The work of Friedrichs assumes that $\langle
\curl \tmmathbf{u} \rangle_{\Om} = 0$, i.e., the curl of the
displacement vector field $\tmmathbf{u}$ in \eqref{eq:K1W1pelementary} is
null-average in $\Om$. Also, it assumes that the region of integration is what
he calls a $\Om$-domain, a more stringent notion than the one of a Lipschitz
domain. In~{\cite{Payne1961}}, Payne and Weinberger investigate the best
possible Korn constants $K_{\Om, 2}$ in balls and show that if Korn's
inequality holds on a finite number of domains, each of which is
$C^2$-diffeomorphic to balls, it also holds on their union.

The extension to domains with the cone property and without the null-average assumption on the
curl of $\tmmathbf{u}$ comes with the work of Gobert~{\cite{Gobert1962}}, who
fruitful uses some previous strategies sketched by Smith to investigate
formally positive integro-differential forms~{\cite{Smith1961}}. 
Later on, Nitsche~{\cite{Nitsche1981}}, who, like Friedrichs, also labels Korn's
original proof as {\emph{doubtful}}, provides an alternative argument for
Korn's second inequality based on the idea of strain-preserving extension
operators which reduces the study of Korn inequalities on bounded Lipschitz
domains to the one on (not necessarily bounded) $C^1$-domains.
All the cited
references are pretty technical and involved. Only with the works of
Kondratiev and Oleinik~{\cite{Kondratiev1989}}, but still, for $p = 2$, a much
shorter, simpler, yet not elementary proof appears.
The relatively recent work of Ciarlet~{\cite{Ciarlet2010}} proposes a new
proof of the second Korn inequality based on a distributional version of the
Saint-Venant lemma. In the same paper, the author also recalls an elegant
proof of the second Korn's inequality from {\cite{Dautray1990}}, yet based on
what he refers to as a {\emph{remarkably difficult-to-prove}} lemma of
J.L.~Lions (for the proof of Lions' lemma, see, e.g., {\cite[Lemma~11.2,
p.~316]{Magenes1958}}, {\cite[Thm~3.2, p.~111]{Duvaut1972}}, and
{\cite{Kesavan_2005}}). 
A more refined version of this lemma is given in \cite[Lemma~2.1]{Dain_2006}, and used to prove a generalized version of Korn's inequality where the linearized strain tensor is replaced by its trace free part. Such an inequality holds, however, only for $n\geqslant 3$.

All the references cited so far cover only the quadratic setting $p = 2$. The
generalization to the general exponent $p \in] 1, \infty [$, i.e.,
inequalities \eqref{eq:K1W1pelementary} and \eqref{eq:K1W1pelementarycc},
introduces new severe difficulties even in the proof of the first Korn's
inequality. The proofs of \eqref{eq:K1W1pelementary} and
\eqref{eq:K1W1pelementarycc} can be found in Mikhlin~{\cite{Mikhlin1965}} and,
under weaker hypotheses, in Mosolov and Mjasnikov~{\cite{Mosolov1971}}.
Both~{\cite{Mikhlin1965}} and {\cite{Mosolov1971}} use Calderón-Zygmund
estimates. A simpler and self-contained proof that does not rely on the theory of singular integrals is given by Kondratiev
and Oleinik~{\cite{Kondrat_ev_1988}}.

Another proof  \cite{Contilecture} of the second Korn's inequality can be given by adapting the arguments used for a nonlinear counterpart thereof, the so-called geometric rigidity inequality \cite{FJM}. One may indeed derive the equivalent statement\footnote{A proof that \eqref{eq:korn-rigid} and  \eqref{eq:K1W1pelementary} are equivalent is given in \cite[Thm~2.3]{Ciarlet2010}: although stated for $p=2$, that proof indeed applies for general $p\neq 1$.} that for all $\tmmathbf{u} \in C^{\infty} (\overline{\Om}, \RR^N)$ there exists a skew-symmetric matrix $A$ such that
\begin{equation}\label{eq:korn-rigid}
\int_{\Om}  | \grad \tmmathbf{u}-A |^p
  \leqslant K_{\Om, p} \int_{\Om}| \gradsym
  \tmmathbf{u} |^p \,.
\end{equation}
The inequality is first proved on a square in the case of {\it harmonic} functions, to which one is reconducted by considering an auxiliary elliptic problem and the related regularity estimates. Then, a fine local-to-global construction making use of a Whitney covering with dyadic squares and of a weighted Poincar\'e inequality (for which some regularity of the boundary is needed) entails the general result.
Let us also mention that if the domain is star-shaped, inequality \eqref{eq:korn-rigid} has been known since the work of Reshetnyak \cite{Reshetnyak_1970}, which is based on a representation of a vector function as the sum of a skew-symmetric affine motion and a potential-type integral only depending on the symmetric gradient. 

We recall that for $p = 1, \infty$ inequalities \eqref{eq:K1W1pelementary} and
\eqref{eq:K1W1pelementarycc} do {\tmem{not}} hold. When $p = 1$, this is a
consequence of a famous result of Ornstein~{\cite{Ornstein_1962}} showing that
a priori estimates for elliptic operators, known to hold in $L^p$ norm ($1 < p
< \infty$), are no longer true in the $L^1$-norm (see also {\cite{Conti2005}}
for a counterexample based on laminates). Similar results of de~Leeuw and
Mirkil~{\cite{Leeuw1962}}, antecedent to the ones of Ornstein, imply that
Korn's inequality cannot be true also if $p = \infty$ (but counterexamples for
the case $p = \infty$ can be easily constructed). From a different
perspective, the result follows from the fact that the so-called $\nabla_2$-
and $\lapl_2$- conditions are sharp requirements for the first and second Korn
inequalities to hold in Orlicz spaces; since the $\Phi$-function associated
with the $L^1$-norm does not satisfy the $\grad_2$-condition, and the
$\Phi$-function associated with the $L^{\infty}$-norm does not satisfy the
$\lapl_2$-condition ~{\cite[Remark~1.3]{Breit2012}}, we get that Korn's
inequality cannot hold when $p = 1, \infty$. 
The results in \cite{Breit2012} have been further extended in \cite{MR3250367} where different, suitably balanced Orlicz norms are considered to cope with the possible failure of either the $\nabla_2$-
or $\lapl_2$- conditions; in \cite{Breit_2017}, where an optimal trace-free inequality is given for $n\geqslant 3$; and in \cite[Prop~4.1]{Conti_2022}, which establishes a general Korn-type inequality for elliptic differential operators. Let us also mention that the even broader class of incompatible Korn inequalities, that is, for vector fields that are not necessarily gradient fields, has been the subject of many recent contributions \cite{Conti_2021, Gmeineder_2021, Gmeineder_2023, Gmeineder_2024, Leoni_2010, MR3294348} that find applications in elastoplasticity models, fluid-mechanical problems, and the analysis of numerical schemes based on the finite element method.

Finally, we recall that Korn domains appear extremely difficult to
characterize. What is known is that not every open set is a Korn domain.
Indeed, as shown by Geymonat and Gilardi in {\cite{Geymonat1998}} and by Weck
in {\cite{Weck1994}}, it is possible to construct domains that satisfy the
segment property, where Korn's second inequality does not hold. 
However, in~\cite{Diening_2010}, it is shown that Korn's inequality holds on a class of domains that is much larger than the class of Lipschitz domains and contains sets that may possess fractal boundaries or internal cusps (external cusps are excluded), the so-called John domains. For finitely connected planar domains, being a John domain is even a necessary condition for a Korn inequality to hold \cite{Jiang2017}.

\subsection{Outline}The rest of the paper is structured as follows.
Section~\ref{sec:setup} presents three Poincar{\'e}-Korn inequalities, valid
even for $p = 1$. One of them is the classical one needed to show that the
symmetric gradient seminorm is a norm in $W^{1, p})$ and
holds for any bounded open subset of $\RR^N$. After that, we state the first
and second Korn's inequality in our geometric setting.

The proofs of the result are given in section~\ref{sec:proofs}. Following the
ideas of {\cite{DiFratta2022}}, we first give elementary proofs of the
Poincar{\'e}-Korn inequalities, whose unique ingredient is the divergence
theorem. Then we focus on our Korn inequalities. Our proofs for the first and
second inequality are essentially the same and based on a lemma
(see~Lemma~\ref{lemmamain}) that we state and prove in the same
section~\ref{sec:proofs}. The lemma allows us to prove the first Korn
inequality on general bounded domains and the second Korn's
inequality on {\tmem{extension domains}} which, in particular, include the class of Lipschitz
domains ---in fact, the class of $(\varepsilon, \delta)$-domains in the sense of Jones (cf.~Remark~\ref{rmk:extedomainE1p}).

\section{Contributions of the present work}\label{sec:setup}
\noindent
The main results of the present work are concise proofs of first and second
Korn's inequality, as well as Poincar{\'e}-Korn inequalities. To state our
results precisely, we need to set up the framework, the mathematical notation,
and the terminology used throughout the paper.

\subsection{Notation}\label{sec:notation}In what follows, we denote by $e_1,
\ldots, e_N$ the standard basis of $\RR^N$ and, for $u \in \RR^N$ we denote by
$| u |$ the euclidean norm, i.e., $| u |^2 \assign (u \cdot e_1)^2 + \cdots +
(u \cdot e_N)^2$. Coherently, the notation $| u |^p$ will stand for the
$p$-power of the Euclidean norm, i.e.,
\begin{equation}
  | u |^p \eqs ((u \cdot e_1)^2 + \cdots + (u \cdot e_N)^2)^{p / 2} .
\end{equation}
For matrices $\Phi, \Psi \in \RR^{N \times N}$, their Frobenius inner product
is defined by $( \Phi \Fsp \Psi ) \assign \tmop{tr} (
\Phi^{\T} \Psi )$, and the associated norm is given by $| \Phi |^2 =
\sum_{i, j = 1}^N (\Phi e_i \cdot e_j)^2$. Coherently, the notation $| \Phi
|^p$ will stand for the $p$-power of $| \Phi |$, i.e.,
\begin{equation}
  | \Phi |^p = \left( \sum_{i, j = 1}^N (\Phi e_i \cdot e_j)^2 \right)^{p / 2}
  .
\end{equation}
We denote by $I$ the identity matrix in $\RR^{N \times N}$. Observe that $| I
| = \sqrt{N}$.

For $\Om \subseteq \RR^N$ open set, we denote by $C^{\infty}_c (\Om, \RR^N)$ 
the space of infinitely differentiable vector fields with
compact support, by $C^{\infty} (\overline{\Om}, \RR^N)$ the
space of infinitely differentiable vector fields that are smooth up to the
boundary, and by $C^{\infty} (\Om, \RR^N)$ the space of
infinitely differentiable vector fields that are smooth {\tmem{in}} $\Om$. We
denote by $\mathcal{D}' (\Om, \RR^N)$ the space of $\RR^N$-valued
distributions on $\Om$. For $\tmmathbf{u} \in \mathcal{D}' ( \Om, \RR^N)$ and $\tmmathbf{\varphi} \in C^{\infty}_c (\Om, \RR^N)$
we denote by $\jdual{\tmmathbf{u}, \tmmathbf{\varphi}}$ the value of
$\tmmathbf{u}$ on $\tmmathbf{\varphi}$.

Similar notation is used when the target space is, e.g., the algebra $\RR^{N
\times N}$ of square matrices of order $N$.

For classical differential operators, the distributional counterparts are the
familiar ones. For $u \in \mathcal{D}' (\Om, \RR)$, $\tmmathbf{u}
\in \mathcal{D}' (\Om, \RR^N)$, $\varphi \in C^{\infty}_c (\Om )$ and $\tmmathbf{\varphi} \in C^{\infty}_c (\Om, \RR^N)$ we set
\begin{equation}
  \jdual{\grad u, \tmmathbf{\varphi}} \assign - \jdual{u, \divv
  \tmmathbf{\varphi}}, \quad \jdual{\divv \tmmathbf{u}, \varphi} \assign -
  \jdual{u, \grad \varphi}, \quad \jdual{\curl \tmmathbf{u},
  \tmmathbf{\varphi}} \assign \jdual{\tmmathbf{u}, \curl \tmmathbf{\varphi}} .
\end{equation}
Also, we set
\begin{equation}
  \jdual{\lapl u, \varphi} \assign \jdual{\divv \grad u, \varphi} \eqs
  \jdual{u, \lapl \varphi}, \quad \jdual{\lapl \tmmathbf{u},
  \tmmathbf{\varphi}} \assign \jdual{\tmmathbf{u}, \lapl \tmmathbf{\varphi}} .
\end{equation}
We need to introduce further differential operators. For $\tmmathbf{u} \in
C^{\infty} (\Om, \RR^N)$, the {\tmem{gradient}} of $\tmmathbf{u}$
is defined by $\grad \tmmathbf{u} \assign ( \nabla \tmmathbf{u} \cdot e_1
\,{\divides}\, \ldots \, {\divides}\, \nabla \tmmathbf{u}
\cdot e_N )$. In other words, $\grad \tmmathbf{u}$ is the matrix in
$\RR^{N \times N}$ whose columns are the gradients of the components of
$\tmmathbf{u}$. Consistently, we denote by $\grad^{\T} \tmmathbf{u} \assign
( \partial_1 \tmmathbf{u} \,{\divides}\, \ldots \,
\,{\divides}\, \partial_N \tmmathbf{u} )$ the Jacobian matrix of
$\tmmathbf{u}$. The {\tmem{symmetric}} and {\tmem{skew-symmetric}} gradient of
$\tmmathbf{u}$ are, respectively, defined by
\begin{equation}
  \gradsym \tmmathbf{u} \assign \frac{1}{2} ( \grad \tmmathbf{u}+
  \grad^{\T} \tmmathbf{u} ), \quad \gradskw \tmmathbf{u} \assign
  \frac{1}{2} ( \grad \tmmathbf{u}- \grad^{\T} \tmmathbf{u} ) .
\end{equation}
Given a matrix-valued field $\tmmathbf{\Phi} \assign (
\tmmathbf{\varphi}_1 \,{\divides}\, \ldots \,{\divides}\,
\tmmathbf{\varphi}_N ) \in C^{\infty} (\Om, \RR^{N \times N})$, we agree that the {\tmem{divergence}} of $\tmmathbf{\Phi}$, denoted
by $\textbf{$\divb$} \tmmathbf{\Phi}$, is the vector field
\begin{equation}
  \textbf{$\divb$} \tmmathbf{\Phi} \assign \sum_{i = 1}^N ( \divv
  \tmmathbf{\varphi}_i ) e_i  \eqs \sum_{j = 1}^N \partial_j
  \tmmathbf{\Phi}^{\T} e_j, \label{eq:exprdivb}
\end{equation}
obtained by applying the classical divergence to the columns of
$\tmmathbf{\Phi}$. With the previous definition, the operator
$\textbf{$\divb$}$ is the formal adjoint of $- \grad$, in the sense that if
$\tmmathbf{u} \in C^{\infty} (\Om, \RR^N)$ and $\tmmathbf{\Phi}
\in C^{\infty}_c (\Om, \RR^{N \times N} )$, then the classical
integration by parts formula holds:
\[ \int_{\Om} \grad \tmmathbf{u} \Fsp \Phi \eqs - \int_{\Om} \tmmathbf{u}
   \cdot \textbf{$\divb$} \tmmathbf{\Phi}. \]
The {\tmem{symmetric}} and {\tmem{skew-symmetric}} divergence of
$\tmmathbf{\Phi}$ are defined by
\begin{align}
  \divbsym \tmmathbf{\Phi} & \assign \frac{1}{2} \divb (
  \tmmathbf{\Phi}+\tmmathbf{\Phi}^{\T} ) \eqs \frac{1}{2} \sum_{j = 1}^N
  \partial_j ( \tmmathbf{\Phi}^{\T} +\tmmathbf{\Phi} ) e_j, \\
  \divbskw \tmmathbf{\Phi} &\assign \frac{1}{2} \divb (
  \tmmathbf{\Phi}-\tmmathbf{\Phi}^{\T} ) \eqs \frac{1}{2} \sum_{j = 1}^N
  \partial_j ( \tmmathbf{\Phi}^{\T} -\tmmathbf{\Phi} ) e_j .
\end{align}
Note that the operators $\textbf{$\divbsym$}$ and $\divbskw$ are the formal
adjoints of $- \gradsym$ and $- \gradskw$. Consistently, if $\tmmathbf{u} \in
\mathcal{D}' (\Om, \RR^N)$, the distributions $\gradsym
\tmmathbf{u} \in \mathcal{D}' (\Om, \RR^{N \times N} )$ and
$\gradskw \tmmathbf{u} \in \mathcal{D}' (\Om, \RR^{N \times N})$
are, respectively, defined for every $\tmmathbf{\Phi} \in C^{\infty}_c (
\Om, \RR^{N \times N} )$ by $\jdual{\gradsym \tmmathbf{u},
\tmmathbf{\Phi}} \assign - \jdual{\tmmathbf{u}, \divbsym \tmmathbf{\Phi}}$ and
$\jdual{\gradskw \tmmathbf{u}, \tmmathbf{\Phi}} \assign - \jdual{\tmmathbf{u},
\divbskw \tmmathbf{\Phi}}$. The definitions of $\gradsym \tmmathbf{u}$ and
$\gradskw \tmmathbf{u}$ are compatible with the smooth setting. If
$\tmmathbf{u} \in C^{\infty} (\Om, \RR^N)$ is a regular
distributions and $\tmmathbf{\Phi} \in C^{\infty}_c (\Om, \RR^{N \times
N})$ then
\begin{equation}
  \jdual{\gradsym \tmmathbf{u}, \tmmathbf{\Phi}} \eqs \int_{\Om} \gradsym
  \tmmathbf{u} \Fsp \tmmathbf{\Phi}, \qquad \jdual{\gradskw \tmmathbf{u},
  \tmmathbf{\Phi}} \eqs \int_{\Om} \gradskw \tmmathbf{u} \Fsp \tmmathbf{\Phi}.
\end{equation}
Similar definitions apply to the distributional version of the operator
$\divbsym$ and $\divbskw$, meaning that if $\tmmathbf{U} \in \mathcal{D}'
(\Om, \RR^{N \times N})$, the distributions $\divbsym
\tmmathbf{U} \in \mathcal{D}' (\Om, \RR^N)$ and $\divbskw
\tmmathbf{U} \in \mathcal{D}' (\Om, \RR^N)$ are defined for every
$\tmmathbf{\varphi} \in C^{\infty}_c (\Om, \RR^N)$ by
$\jdual{\divbsym \tmmathbf{U}, \tmmathbf{\varphi}} \assign -
\jdual{\tmmathbf{U}, \gradsym \tmmathbf{\Phi}}$ and $\jdual{\divbskw
\tmmathbf{U}, \tmmathbf{\varphi}} \assign - \jdual{\tmmathbf{U}, \gradskw
\tmmathbf{\Phi}}$. Simple computations show that for $\tmmathbf{u} \in
C^{\infty} (\Om, \RR^N)$ there holds
\begin{align}
  - \lapl \tmmathbf{u} & \eqs  - \grad \divv \tmmathbf{u}- 2 \divbskw
  \gradskw \tmmathbf{u},  \label{eq:HodgeN}\\
  - \lapl \tmmathbf{u} & \eqs  + \grad \divv \tmmathbf{u}- 2 \divbsym
  \gradsym \tmmathbf{u}, 
\end{align}
and very same relations hold if $\tmmathbf{u} \in \mathcal{D}' (\Om,
\RR^N )$. Note that if $N = 3$, we have $- 2 \divbskw \gradskw
\tmmathbf{u}= \curl  \curl \tmmathbf{u}$ and \eqref{eq:HodgeN} reduces to the
classical Helmholtz decomposition in $\RR^3$.

\subsection{Main results}Our first result concerns Poincar{\'e}-Korn
inequalities in bounded and unbounded domains. These inequalities are
essential tools in showing the coercivity of variational problems of
elasticity. One of them ({\tmabbr{cf.}}~\eqref{eq:PKCc3}) is a propedeutic
tool for the proof of the first Korn's inequality: it shows that the symmetric
gradient seminorm is actually a norm on $C^{\infty}_c (\Om, \RR^N)$.

\begin{theorem}[{\tmname{Poincar{\'e}-Korn inequalities}}]
  \label{thm:mainPK}Let $\Om \subseteq \RR^N$ be an open set and $1 \leqslant
  p < \infty$. The following Poincar{\'e}-Korn inequalities hold.
  \begin{enumerateroman}
    \item There exists a constant $C_{p, N} > 0$ depending only on $p$ and
    $N$, such that
    \begin{equation}
      \int_{\Om} | \tmmathbf{u} |^p \leqslant C_{p, N}^p \int_{\Om} | x |^p 
      | \gradsym \tmmathbf{u} |^p \qquad \forall \tmmathbf{u} \in
      C^{\infty}_c (\Om, \RR^N) . \label{eq:PKCc1}
    \end{equation}
    In particular, if $\Om$ is also bounded, then there exists a constant
    $\kappa_{p, \Om} > 0$ depending only on $p$ and $\Om$, such that
    \begin{equation}
      \int_{\Om} | \tmmathbf{u} |^p \leqslant \kappa_{p, \Om}^p \int_{\Om}
      | \gradsym \tmmathbf{u} |^p \qquad \forall \tmmathbf{u} \in
      C^{\infty}_c (\Om, \RR^N) . \label{eq:PKCc3}
    \end{equation}
    \item If $\Om$ is a bounded Lipschitz domain, then there exist constants
    $\kappa_{p, \Om}, \kappa_{p, \partial \Om} > 0$, depending only on $p$ and
    $\Om$ such that
    \begin{equation}
      \int_{\Om} | \tmmathbf{u} |^p \leqslant \kappa_{p, \Om}^p \int_{\Om}
      | \gradsym \tmmathbf{u} |^p + \kappa_{p, \partial \Om}
      \int_{\partial \Om} | \tmmathbf{u} |^p \qquad \forall \tmmathbf{u} \in
      C^{\infty} (\overline{\Om}, \RR^N) . \label{eq:PKCc4}
    \end{equation}
  \end{enumerateroman}
\end{theorem}

We give an elementary proof of Theorem~\ref{thm:mainPK} in
section~\ref{sec:PoincKornproofs}. Here, instead, we make some remarks. First,
we stress that in contrast to Korn's inequality which does not hold for $p =
1$, the Poincar{\'e}-Korn inequality holds even for $p = 1$ and
{\tmem{arbitrary}} open sets. Second, we note that {\tmem{ii.}} implies
\eqref{eq:PKCc3}. Indeed, if $\Om$ is bounded but {\tmem{not}} Lipschitz, one
can extend by zero an element $\tmmathbf{u} \in C_c^{\infty} (\Om, \RR^N
)$ to an element of $\tmmathbf{u}_0 \in C_c^{\infty} ( B, \RR^N
)$ where $B$ is a ball containing $\Om$, and evaluate \eqref{eq:PKCc4}
on $\tmmathbf{u}_0$ (and with $\Om$ replaced by $B$) to get back to
\eqref{eq:PKCc3}. We opted to present \eqref{eq:PKCc3} as a particular
case of \eqref{eq:PKCc1}. 
Finally, it will be apparent from
the proof that one can take
\begin{equation}
  C_{p, N} \assign (2+ | p - 2 | + \sqrt{N}) \cdot \frac{p}{p +
  N}, \quad \kappa_{p, \Om} \assign ( \tmop{diam} \Om ) C_{p, N},
  \quad \kappa_{p, \partial \Om} \assign p \left( \frac{p + 1}{p + N} \right)
  ( \tmop{diam} \Om ) . \label{eq:PKCc2}
\end{equation}
However, there is no pretense of optimality in the previous constants.

Our second contribution is the main objective of our paper and concerns a
concise point of view on the first and second Korn inequality. To state the
results in the generality we want, we need to introduce the proper functional
setting.

For every $1 \leqslant p < \infty$, we denote by $E^{1, p} (\Om, \RR^N)$ the vector space of $L^p$-vector-fields, whose symmetric gradient, in
the distributional sense, is still in $L^p$, i.e.,
\begin{equation}
  E^{1, p} (\Om, \RR^N) \assign \{ \tmmathbf{u} \in L^p
  (\Om, \RR^N) \of \gradsym \tmmathbf{u} \in L^p (\Om,
  \RR^{N \times N}) \} .
\end{equation}
It is standard to show (see, e.g., {\cite{Temam2019}}) that $E^{1, p} (
\Om )$ is a Banach space when endowed with the norm
\begin{equation}
  \| \tmmathbf{u} \|_{E^{1, p}} \assign ( \| \tmmathbf{u} \|^p_{L^p} +
  \| \gradsym \tmmathbf{u} \|^p_{L^p} )^{1 / p} = \left(
  \int_{\Om} | \tmmathbf{u} |^p + | \gradsym \tmmathbf{u} |^p
  \right)^{1 / p} . \label{eq:normE1p}
\end{equation}
Also, the Poincar{\'e}-Korn inequality \eqref{eq:PKCc3} assures that on the
Banach subspace of $E^{1, p} (\Om, \RR^N)$, given by
\begin{equation}
  E^{1, p}_0 (\Om, \RR^N) \assign \overline{C_c^{\infty} (
  \Om, \RR^N )},
\end{equation}
i.e., defined as the closure of $C_c^{\infty} (\Om, \RR^N)$ in
$E^{1, p} (\Om, \RR^N)$, the seminorm $\| \gradsym
\tmmathbf{u} \|_{L^p}$ is actually a norm equivalent to
\eqref{eq:normE1p}. It is standard to show that, as happens for the space
$W^{1, p}_0 (\Om, \RR^N)$, even $E^{1, p}_0 (\Om, \RR^N
)$ has the extension-by-zero property. Namely, if $\tmmathbf{u} \in
E^{1, p}_0 (\Om, \RR^N)$ then
\begin{equation}
  \tilde{\tmmathbf{u}} \in E^{1, p} ( \RR^N, \RR^N ) \quad
  \text{and} \quad \gradsym \tilde{\tmmathbf{u}} = \widetilde{\gradsym
  \tmmathbf{u}},
\end{equation}
where we denoted by $\tilde{\tmmathbf{u}}$ and $\widetilde{\gradsym
\tmmathbf{u}}$ the extensions by zero of $\tmmathbf{u}$ and $\gradsym
\tmmathbf{u}$ to the whole of $\RR^N$. Indeed, if $\tmmathbf{u} \in E^{1, p}_0
(\Om, \RR^N)$ and $(\tmmathbf{\varphi}_n)_{n \in \NN}$ is a
sequence in $C_c^{\infty} (\Om, \RR^N)$ such that
$\tmmathbf{\varphi}_n \to \tmmathbf{u}$ in $E^{1, p}_0 (\Om,
\RR^N )$, then for every $\tmmathbf{\Phi} \in C^{\infty}_c ( \RR^N,
\RR^{N \times N})$ we have that
\begin{align}
  \jdual{\gradsym \tilde{\tmmathbf{u}}, \tmmathbf{\Phi}}_{\RR^N} & \eqs  -
  \int_{\Om} \tmmathbf{u} \cdot \divbsym \tmmathbf{\Phi}  \eqs -
  \lim_{n \to \infty} \int_{\Om} \tmmathbf{\varphi}_n \cdot \divbsym
  \tmmathbf{\Phi} \; \; \eqs \; \; \lim_{n \to \infty} \int_{\Om}
  \grad_{\tmop{sym}} \tmmathbf{\varphi}_n \of \tmmathbf{\Phi} \nonumber\\
  & \eqs  \int_{\RR^N} \widetilde{\gradsym \tmmathbf{u}} \of \tmmathbf{\Phi}
   \eqs  \jdual{\widetilde{\gradsym \tmmathbf{u}},
  \tmmathbf{\Phi}}_{\RR^N} .  \label{eq:extbyzero}
\end{align}
But the extension-by-zero property is not the only feature that the spaces
$E^{1, p}_0 (\Om, \RR^N)$ and $W^{1, p}_0 (\Om, \RR^N)$ share. 
Indeed, as we show in Theorem~\ref{thm:mainK1} below, when $1
< p < \infty$ and $\Om$ is a bounded open set, then  $E^{1, p}_0 (\Om, \RR^N) =
W^{1, p}_0 (\Om, \RR^N)$ and their associated norms are
equivalent. Also,  one has $E^{1, p} (\Om, \RR^N) = W^{1, p} (\Om,\RR^N)$ when $1 < p < \infty$ and $\Om$ is what we call an {\tmem{extension
domain}}.

\begin{definition}
  \label{def:restext}We say that a bounded open set $\Om$ is an extension
  domain {\opt}for $E^{1, p}${\cpt} when any element $\tmmathbf{u} \in E^{1,
  p} (\Om, \RR^N)$ can be extended to an element of $E^{1, p}
  ( \RR^N, \RR^N )$. 
\end{definition}

\begin{remark} \label{rmk:extedomainE1p}
Generally speaking, $E^{1, p}$-extension domains can be
  investigated by transposing the known techniques from $W^{1, p}$ to $E^{1,
  p}$. For example, the technique in \cite{Nitsche1981}, which suitably adapts the standard extension by reflection to the setting of symmetrized gradients, shows that Lipschitz domains are $E^{1, p}$-extension domains. Recently, it has been shown that $(\varepsilon, \delta)$-domains in the sense of Jones are $E^{1, p}$-extension domains. For that, we refer the reader to \cite[Thm~6.2]{Breit_2021} which generalizes previous results contained in \cite{Gmeineder_2019}. Finally, we observe that a notion of $E^{1, p}_0$-extension domain would not be
  interesting because, as we have shown in \eqref{eq:extbyzero}, every bounded
  open set would be of this type.
\end{remark}

\begin{theorem}[{\tmname{first and second Korn's inequality}}]
  \label{thm:mainK1} For $1 < p< \infty$, the following assertions hold:
  \begin{enumerateroman}
    \item If $\Om \subseteq \RR^N$ is an extension domain, then $E^{1, p} (\Om, \RR^N
    ) = W^{1, p} (\Om, \RR^N)$ and there exists a constant
    $K_{\Om, p} > 0$ such that
    \begin{equation}
      \int_{\Om} | \tmmathbf{u} |^p + | \grad \tmmathbf{u} |^p
      \leqslant K_{\Om, p} \int_{\Om} | \tmmathbf{u} |^p + | \gradsym
      \tmmathbf{u} |^p \quad \forall \tmmathbf{u} \in W^{1, p} (
      \Om, \RR^N ) . \label{eq:K1W1p}
    \end{equation}
    \item If $\Om \subseteq \RR^N$ is an open set, then $E^{1, p}_0 ( \Om
    , \RR^N) = W^{1, p}_0 ( \Om , \RR^N)$ and there exists a constant
    $K_{\Om, p} > 0$ such that
    \begin{equation}
      \int_{\Om} | \grad \tmmathbf{u} |^p \leqslant K_{\Om, p}
      \int_{\Om} | \gradsym \tmmathbf{u} |^p \quad \forall
      \tmmathbf{u} \in W^{1, p}_0 (\Om, \RR^N) .
      \label{eq:K1W1p0}
    \end{equation}
  \end{enumerateroman}
\end{theorem}

\begin{remark}
  \label{rmk:forp1inf} As already remarked in the introduction, as long as one is not interested in optimal constants, a simple scaling argument on balls gives that the constant $K_{\Om, p}$ in equation \eqref{eq:K1W1p0} can always be taken to be independent of $\Om$.  
  Also, as already recalled in section~\ref{sec:stateoftheart},
  Korn inequalities are no longer true when $p = 1$ and the space $W^{1, 1}
  (\Om, \RR^N)$ and $E^{1, 1} (\Om, \RR^N)$ are
  different. The space $E^{1, 1} (\Om, \RR^N)$ is often denoted
  by $\tmop{LD} (\Om, \RR^N)$ and plays a pivotal role in the
  mathematical theory of plasticity. We refer the reader to {\cite[Chap.
  II]{Temam2019}} for some of the main properties of the space $\tmop{LD}
  (\Om, \RR^N)$.
\end{remark}

The concise proof of Theorem~\ref{thm:mainK1} is given at the beginning of
section~\ref{sec:proofs} and relies on the following lemma, whose proof is
deferred to the end of section~\ref{sec:proofs}.

\begin{lemma}
  \label{lemmamain}Let $\Om$ be an open set and $1 < p < \infty$. If
  $\tmmathbf{u} \in E^{1, p} (\Om, \RR^N)$, then $\tmmathbf{u}
  \in W_{\tmop{loc}}^{1, p} (\Om, \RR^N)$.
\end{lemma}

The approach we follow to prove the first and second Korn's inequality
appeared concise and effective in showing the regularity of distributional
solutions of the Poisson equation~{\cite{Di_Fratta_2020bb}}. Part of our
arguments has a nonempty intersection with the one Duvaut and Lions gave in
{\cite[Thm~3.2, p.~111]{Duvaut1972}}. However, our approach avoids the
{\tmem{remarkably difficult-to-prove}} lemma of
J.L.~Lions~{\cite{Dautray1990}} and, for $p = 2$, one only needs the classical
Riesz representation theorem in Hilbert space and the Weyl lemma on the
$C^{\infty}$-regularity of harmonic distribution (for which short and
elementary proofs are available~{\cite{Beltrami_1968}}). Even for general $1 <
p < \infty$, we rely on the same ingredients: the $q$-Riesz representation
theorem~{\cite[pp.~10-11]{Simader1996}}, with $q$ such that $1 / q + 1 / p =
1$, and Weyl lemma. We are unaware of a simple proof of the $q$-Riesz
representation theorem when $q \neq 2$. However, our proof requires invoking
the result only on balls.

\section{First and second Korn's inequality: proofs of
Theorem~\ref{thm:mainK1} and Lemma \ref{lemmamain} }\label{sec:proofs}
\noindent
\subsection{Proof of Theorem~\ref{thm:mainK1}.{\tmem{i}}}  We first show that if
$1 < p < \infty$ then $W^{1, p} (\Om, \RR^N) = E^{1, p}(
\Om, \RR^N )$. By the very definition, we have $W^{1, p} ( \Om
) \subseteq E^{1, p} (\Om, \RR^N)$, even for $p=1$ and general open sets. Vice versa, if $1<p<\infty$, $\Om$ is
an extension domain, and $\tmmathbf{u}
\in E^{1, p} (\Om, \RR^N)$, then its extension
$\tilde{\tmmathbf{u}} \in E^{1, p} ( \RR^N, \RR^N )$ is, by
Lemma~\ref{lemmamain}, in $W^{1, p}_{\tmop{loc}} ( \RR^N, \RR^N )$.
Therefore, $\tmmathbf{u} \in W^{1, p} (\Om, \RR^N)$.

Next, we observe that the inclusion map
\begin{equation}
  \jmath \of \tmmathbf{u} \in W^{1, p} (\Om, \RR^N) \mapsto
  \tmmathbf{u} \in E^{1, p} (\Om, \RR^N)  \label{eq:jembedd}
\end{equation}
is continuous for any $1 \leqslant p < \infty$, because of the inequality
$\| \gradsym \tmmathbf{u} \|_{L^p} \leqslant \| \grad
\tmmathbf{u} \|_{L^p}$, and, if $p > 1$, also a
surjection of $W^{1, p} (\Om, \RR^N)$ onto $E^{1, p} ( \Om
)$. Therefore, by the open mapping theorem, $\jmath$ is a
topological isomorphism of $W^{1, p} (\Om, \RR^N)$ onto $E^{1, p}
(\Om, \RR^N)$, and this implies \eqref{eq:K1W1p}.
{\jqed}

\subsection{Proof of Theorem~\ref{thm:mainK1}.{\tmem{ii}}}It is sufficient to
show that if $1 < p < \infty$ then $W_0^{1, p} (\Om, \RR^N) = E^{1, p}_0 (\Om,
\RR^N)$. Indeed, after that, inequality \eqref{eq:K1W1p0} follows from the
fact that the inclusion map $\jmath \of \tmmathbf{u} \in W^{1, p}_0 (\Om,
\RR^N) \mapsto \tmmathbf{u} \in E^{1, p}_0 (\Om, \RR^N)$ is a continuous
surjection of the Banach space $(W^{1, p}_0 (\Om, \RR^N), \| \grad \cdot \|_{L^p})$ onto the Banach space $(E^{1, p}_0
(\Om, \RR^N), \| \gradsym \cdot \|_{L^p})$.

The inclusion $W_0^{1, p} (\Om, \RR^N) \subseteq E^{1, p}_0 (\Om, \RR^N)$ is
trivial and holds even for $p = 1$. Indeed, with $\jmath$ given by
\eqref{eq:jembedd}, we get that
\[ W^{1, p}_0 (\Om, \RR^N) = \jmath (\overline{C_c^{\infty} (\Om, \RR^N)}
   \subseteq \overline{\jmath (C_c^{\infty} (\Om, \RR^N))} = E^{1, p}_0 (\Om,
   \RR^N), \]
where the first closure is meant with respect to the topology induced by the
$W^{1, p}$-norm, while the second refers to the topology induced by
the $E^{1, p}$-norm.
For the inclusion $E^{1, p}_0 (\Om, \RR^N) \subseteq W_0^{1, p} (\Om, \RR^N)$
we proceed in two steps. In \tmtextsc{Step}\tmrsub{$1$}, we assume that $\Om$
is bounded and of class $C^1$, and then, in \tmtextsc{Step}\tmrsub{$2$}, we show how to
remove this hypothesis.
\smallskip

\noindent\tmtextsc{Step}\tmrsub{$1$}. If $\tmmathbf{u} \in E^{1, p}_0 (\Om,
\RR^N)$ then its extension by zero to the whole of $\RR^N$, let us denote it
by $\tilde{\tmmathbf{u}}$, is in $E^{1, p} (\RR^N, \RR^N)$ and has
distributional support in $\Om$. Hence, by Lemma~\ref{lemmamain}, we have that
\begin{equation}
  \tilde{\tmmathbf{u}} \in W^{1, p}_{\tmop{loc}} (\RR^N, \RR^N), \quad
  \tmop{supp} \tilde{\tmmathbf{u}} \subseteq \overline{\Om} .
\end{equation}
Since $\Om$ is a domain of class $C^1$, we conclude\footnote{The motivation
for a first step under the assumption that $\Om$ is of class $C^1$ comes from
the fact that such domains satisfy the property that if $\tilde{\tmmathbf{u}}
\in W^{1, p} (\RR^N, \RR^N)$ and $\tmop{supp} \tilde{\tmmathbf{u}} \subseteq
\overline{\Om}$, then $\tilde{\tmmathbf{u}} \in W_0^{1, p} (\Om, \RR^N)$. The
result is proved, e.g., in {\cite[Thm.~2, p.~273]{Evans2010}}. It is also
known that domains of class $C^0$ still assure the same result (see
{\cite[sec.~2.1]{Davoli2022a}} and {\cite{Chandler_Wilde_2017}} for details).
However, to the best of our knowledge, it is unknown how general this class
can be, and that is why we have to rely on a second step. The idea of arguing
in two steps to have a result valid for a general open set was suggested to us
by the anonymous referee. We thank the anonymous referee for this and other
valuable suggestions.} that $\tmmathbf{u} \in W^{1, p}_0 (\Om, \RR^N)$.

\smallskip

\noindent
\tmtextsc{Step}\tmrsub{$2$}. Let $\Om$ be a bounded open set,
$\tmmathbf{u} \in E^{1, p}_0 (\Om, \RR^N)$, and $(\tmmathbf{u}_n)_{n \in
\mathbb{N}}$ a sequence in
$C_c^{\infty} (\Om, \RR^N)$ converging to $\tmmathbf{u}$ in $E^{1, p}_0 (\Om,
\RR^N)$. It is sufficient to prove that $(\tmmathbf{u}_n)_{n \in
\mathbb{N}}$ is a Cauchy sequence in
$W^{1, p}_0 (\Om, \RR^N)$ because, after that, we have \ $\tmmathbf{u}_n
\rightarrow \tmmathbf{v}$ in $L^p$ for some $\tmmathbf{v} \in W^{1, p}_0 (\Om,
\RR^N)$ and, therefore, $\tmmathbf{u}=\tmmathbf{v}$. To that end, let $\tilde{\tmmathbf{u}}_n \in C^{\infty}_c ( B, \RR^N )$ denote
the trivial extensions by zero of $\tmmathbf{u}_n$ to a ball $B$ that includes
$\Om$. By \tmtextsc{Step}\tmrsub{$1$}, Korn's inequality
\eqref{eq:K1W1p0} holds for domains of class $C^1$ (in particular, for $B$). Therefore, for every $n,m\in\NN,n>m$ there holds
\begin{align}
  \lim_{n \rightarrow \infty} \sup_{m > n} \| \tmmathbf{u}_n -\tmmathbf{u}_m
  \|_{W^{1, p}_0 (\Omega)} & =  \lim_{n \rightarrow \infty} \sup_{m > n} \|
  \tilde{\tmmathbf{u}}_n - \tilde{\tmmathbf{u}}_m \|_{W^{1, p}_0 (B)} \nonumber\\
   & \leqslant  K_{B, p} \cdot \lim_{n \rightarrow \infty} \sup_{m > n} \|
  \tilde{\tmmathbf{u}}_n - \tilde{\tmmathbf{u}}_m \|_{E^{1, p}_0 (B)}  \nonumber\\
   & =  K_{B, p} \cdot \lim_{n \rightarrow \infty} \sup_{m > n} \|
  \tmmathbf{u}_n -\tmmathbf{u}_m \|_{E^{1, p}_0 (\Omega)} 
  =  0. \nonumber
\end{align}
This completes the proof if the case where $\Om$ is bounded. The case of $\Om$ being unbounded follows from the scale invariance of Korn's constant on balls.{\jqed}

\subsection{Proof of Lemma~\ref{lemmamain}}The argument we present is inspired
by the one in~{\cite{Di_Fratta_2020bb}}. We first give the proof for $N = 3$,
and show afterward how to modify the proof for general $N$. The proof for
general $N$ is as easy as in the case $N = 3$, but for $N = 3$, we have the
familiar $\curl$ operator at our disposal, and we don't need further
notation.{\medskip}

{\noindent}{\tmstrong{Proof of Lemma~\ref{lemmamain} for $N = 3$.}} Let
$\tmmathbf{u} \in E^{1, p} (\Om, \RR^N)$ and consider a ball $B
\subseteq \Om$. By the $q$-exponent version of the Riesz representation
theorem {\cite[pp.~10-11]{Simader1996}} there exists an element $\tmmathbf{v}
\in W^{1, p}_0 ( B, \RR^N )$, $p > 1$, such that
\begin{equation}
  - \frac{1}{2} \jdual{\lapl \tmmathbf{v}, \tmmathbf{\varphi}} \eqs
  \jdual{\gradsym \tmmathbf{u}, \gradsym \tmmathbf{\varphi}} \quad \forall
  \tmmathbf{\varphi} \in C_c^{\infty} ( B, \RR^N ) .
  \label{eq:rrt1loc}
\end{equation}
Indeed, the right-hand side of \eqref{eq:rrt1loc} defines a bounded linear
functional on $W^{1, q}_0( B, \RR^N )$ ---note that for $q = p =
2$, we only need the classical Riesz representation theorem in Hilbert spaces.

Also, simple algebra and integration by parts formula, give that $2
\jdual{\gradsym \tmmathbf{u}, \gradsym \tmmathbf{\varphi}} = - \jdual{\lapl
\tmmathbf{u}+ \grad \divv \tmmathbf{u}, \tmmathbf{\varphi}}$ for every
$\tmmathbf{\varphi} \in C_c^{\infty} ( B, \RR^N )$. Therefore,
by combining this equality with \eqref{eq:rrt1loc}, we get that there exists
$\tmmathbf{v} \in W^{1, p}_0 ( B, \RR^N)$, $p > 1$, such that
\begin{equation}
  - \lapl \tmmathbf{v}= - ( \lapl \tmmathbf{u}+ \grad \divv \tmmathbf{u}
  ) \quad \text{in } \mathcal{D}'( B, \RR^N ) .
  \label{eq:repcurlcurlloc}
\end{equation}
In particular, applying the divergence operator to both sides of
\eqref{eq:repcurlcurlloc} we infer that the distribution $\divv
(\tmmathbf{u}-\tmmathbf{v}/ 2)$ is harmonic in $B$, and this, by Weyl's
lemma~{\cite{Beltrami_1968}}, assures that
\begin{equation}
  \divv (\tmmathbf{u}-\tmmathbf{v}/ 2) \in C^{\infty} (B) .
  \label{eq:divCinfty}
\end{equation}
Similarly, applying the $\curl$ operator to both sides of
\eqref{eq:repcurlcurlloc} we deduce that the distribution $\curl 
(\tmmathbf{u}-\tmmathbf{v})$ is harmonic in $B$ and, therefore, that
\begin{equation}
  \curl  (\tmmathbf{u}-\tmmathbf{v}) \in C^{\infty} ( B, \RR^N ) .
\end{equation}
But then
\begin{align}
  - \lapl (\tmmathbf{u}-\tmmathbf{v}) & \eqs  \curl  \curl 
  (\tmmathbf{u}-\tmmathbf{v}) - \grad \divv (\tmmathbf{u}-\tmmathbf{v}/ 2) +
  \grad \divv (\tmmathbf{v}/ 2) \nonumber\\
  & \in  C^{\infty} ( B, \RR^N ) + C^{\infty} ( B, \RR^N
  ) + \grad \divv ( W^{1, p}_0 ( B, \RR^N ) )
  \nonumber\\
  & \in  C^{\infty} ( B, \RR^N ) + C^{\infty}( B, \RR^N
  ) + W^{- 1, q'} ( B, \RR^N ) 
\end{align}
and, therefore, given that $- \mathLaplace \tmmathbf{v} \in W^{- 1, q'} (
B, \RR^N )$ we obtain that
\begin{equation}
  - \lapl \tmmathbf{u} \in C^{\infty} ( B, \RR^N ) + W^{- 1, q'}
  ( B, \RR^N ) . \label{eq:reglaplu}
\end{equation}
The previous relation completes the proof. Indeed, shrinking eventually the ball $B$, we get that $- \lapl \tmmathbf{u} \in C^{\infty} ( \bar{B}, \RR^N ) + W^{- 1, q'} (B, \RR^N ) = W^{- 1, q'} ( B, \RR^N )$. Applying again the $q$-exponent version of the Riesz representation theorem, we deduce the existence of a vector field $\tmmathbf{w} \in W^{1, p}_0 ( B, \RR^N )$ such that $- \lapl (\tmmathbf{u}-\tmmathbf{w}) = 0$ in $\mathcal{D}' (B,\RR^N )$. But then, as before, harmonicity implies that $\tmmathbf{u} \in C^{\infty} ( B, \RR^N ) + W^{1, p}_0 ( B, \RR^N ) \subseteq W^{1, p}_{\tmop{loc}} ( B, \RR^N )$ and we conclude.
{\smallskip}

{\noindent}{\tmstrong{Proof of Lemma~\ref{lemmamain} for general $N$.}} We can
resume from \eqref{eq:divCinfty}. Applying the $\gradskw$ operator to both
sides of \eqref{eq:repcurlcurlloc}, we deduce that the distribution $\gradskw
(\tmmathbf{v}-\tmmathbf{u})$ is harmonic in $B$ and, therefore, that $\gradskw
(\tmmathbf{v}-\tmmathbf{u})$ is smooth. But then
\begin{align}
  - \lapl (\tmmathbf{u}-\tmmathbf{v}) & \eqs  - 2 \divbskw \gradskw
  (\tmmathbf{u}-\tmmathbf{v}) - \grad \divv (\tmmathbf{u}-\tmmathbf{v}/ 2) +
  \grad \divv (\tmmathbf{v}/ 2) \nonumber\\
  & \in  C^{\infty} ( B, \RR^N ) + C^{\infty} ( B, \RR^N
  ) + \grad \divv ( W^{1, p}_0 ( B, \RR^N ) )
  \nonumber\\
  & \in  C^{\infty} ( B, \RR^N ) + C^{\infty} ( B, \RR^N
  ) + W^{- 1, q'} ( B, \RR^N ) 
\end{align}
and, therefore, given that $- \mathLaplace \tmmathbf{v} \in W^{- 1, q'} (
B, \RR^N )$, we are back to \eqref{eq:reglaplu}, and from that point on,
the proof is the same as for $N = 3$. {\jqed}

\section{Proof of Theorem~\ref{thm:mainPK}: Poincar{\'e}-Korn
inequalities}\label{sec:PoincKornproofs}
\noindent
Following the ideas of {\cite{DiFratta2022}}, we first give elementary proofs
of the Poincar{\'e}-Korn inequalities, whose unique ingredient is the
divergence theorem.

{\smallskip}

{\noindent}{\tmstrong{Proof of Theorem~\ref{thm:mainPK}}.} For
every $\varepsilon > 0$, we consider the vector field
$\tmmathbf{u}_{\varepsilon} : \Om \to \RR^{N + 1}$ defined by
$\tmmathbf{u}_{\varepsilon} \assign (\tmmathbf{u}, \varepsilon)$. Note that $|
\tmmathbf{u}_{\varepsilon} | = (| \tmmathbf{u} |^2 + \varepsilon^2)^{1 / 2}$
and, therefore, $| \tmmathbf{u}_{\varepsilon} |^{\alpha} \in C^{\infty} (
\Om, \RR^N )$ for every $\alpha \in \RR$. Here, with a convenient abuse
of notation, we used the same symbol to denote the Euclidean norms in $\RR^{N
+ 1}$ and $\RR^N$. It is easy to check that the following equalities hold
\begin{align}
  | \tmmathbf{u}_{\varepsilon} |^{p - 2} \tmmathbf{u} \cdot ( \grad^{\T}
  \tmmathbf{u}x ) & \eqs  \frac{1}{p} \divv (|
  \tmmathbf{u}_{\varepsilon} |^p x) - \frac{N}{p} | \tmmathbf{u}_{\varepsilon}
  |^p, \\
  | \tmmathbf{u}_{\varepsilon} |^{p - 2} \tmmathbf{u} \cdot ( \grad
  \tmmathbf{u}x ) & \eqs  \divv ((\tmmathbf{u} \cdot x) |
  \tmmathbf{u}_{\varepsilon} |^{p - 2} \tmmathbf{u}) - (\tmmathbf{u} \cdot x)
  \divv (| \tmmathbf{u}_{\varepsilon} |^{p - 2} \tmmathbf{u}) - | \tmmathbf{u}
  |^2 | \tmmathbf{u}_{\varepsilon} |^{p - 2}, 
\end{align}
from which, summing term by term, we get
\begin{align}
  \left( \frac{N}{p} + \frac{| \tmmathbf{u} |^2}{| \tmmathbf{u}_{\varepsilon}
  |^2} \right) | \tmmathbf{u}_{\varepsilon} |^p  & \eqs  - 2  |
  \tmmathbf{u}_{\varepsilon} |^{p - 1}  \frac{\tmmathbf{u}}{|
  \tmmathbf{u}_{\varepsilon} |} \cdot ( \gradsym \tmmathbf{u}\, x ) -
  (\tmmathbf{u} \cdot x) \divv (| \tmmathbf{u}_{\varepsilon} |^{p - 2}
  \tmmathbf{u}) \nonumber\\
  &   \quad \quad \quad  \quad \quad \quad + \frac{1}{p} \divv (|
  \tmmathbf{u}_{\varepsilon} |^p x) + \divv \left ( |
  \tmmathbf{u}_{\varepsilon} |^p \left( \frac{\tmmathbf{u}}{|
  \tmmathbf{u}_{\varepsilon} |} \cdot x \right) \frac{\tmmathbf{u}}{|
  \tmmathbf{u}_{\varepsilon} |} \right) .  \label{eq:fundrel}
\end{align}
The previous relation is the key ingredient in the proof of
Theorem~\ref{thm:mainPK}.

{\smallskip}

{\noindent}{{\tmstrong{Proof of
Theorem~\ref{thm:mainPK}.{\tmem{i.}}}}} {\sc{Proof of} \eqref{eq:PKCc1}.}
If $\tmmathbf{u} \in C^{\infty}_c (\Om, \RR^N)$, integrating on
$\Om$ both sides \eqref{eq:fundrel}, and then passing to the limit for
$\varepsilon \to 0$, we get that
\begin{align}
  \left( \frac{p + N}{p} \right) \int_{\Om} | \tmmathbf{u} |^p & \leqslant 
  2 \int_{\Om} | \tmmathbf{u} |^{p - 1} | \gradsym \tmmathbf{u} |  |
  x | + \liminf_{\varepsilon \to 0} \int_{\Om} | \tmmathbf{u} | 
  | \divv (| \tmmathbf{u}_{\varepsilon} |^{p - 2} \tmmathbf{u}) | 
  | x |,  \label{eq:temptoest1}
\end{align}
and, we are left to evaluate the last term on the right-hand side of the
previous relation \eqref{eq:temptoest1}. For that, we observe that
\begin{equation}
  \divv (| \tmmathbf{u}_{\varepsilon} |^{p - 2} \tmmathbf{u})  \eqs  |
  \tmmathbf{u}_{\varepsilon} |^{p - 2} \left[ (p - 2)  ( \gradsym
  \tmmathbf{u} ) \frac{\tmmathbf{u}}{| \tmmathbf{u}_{\varepsilon} |}
  \cdot \frac{\tmmathbf{u}}{| \tmmathbf{u}_{\varepsilon} |} + \divv
  \tmmathbf{u} \right] 
\end{equation}
and, therefore, given that $| \divv \tmmathbf{u} | \eqs |
\gradsym \tmmathbf{u} \of I | \leqslant \sqrt{N} | \gradsym
\tmmathbf{u} |$, we have the estimate
\begin{equation}
  \liminf_{\varepsilon \to 0} \int_{\Om} | \tmmathbf{u} |  |
  \divv (| \tmmathbf{u}_{\varepsilon} |^{p - 2} \tmmathbf{u}) |  | x |
  \leqslant ( | p - 2 | + \sqrt{N} ) \int_{\Om} | \tmmathbf{u} |^{p
  - 1} | \gradsym \tmmathbf{u} | . \label{eq:temptoest2}
\end{equation}
Estimating \eqref{eq:temptoest1} with \eqref{eq:temptoest2} we get that
\[ \left( \frac{N + p}{p} \right) \int_{\Om} | \tmmathbf{u} |^p \leqslant
   (2 + | p - 2 | + \sqrt{N} ) \int_{\Om} | \tmmathbf{u} |^{p - 1}
   | \gradsym \tmmathbf{u} |  | x | \, ,\]
which establishes \eqref{eq:PKCc1} when $p = 1$. For $p > 1$, H{\"o}lder
inequality gives the Poincar{\'e}-Korn inequality under the weighted form
\eqref{eq:PKCc1}.

{\noindent}{\tmname{Proof of \eqref{eq:PKCc3}.}} It is sufficient to note that
the proof of \eqref{eq:PKCc1} works verbatim if one replaces everywhere $x$
with $x - x_0$, $x_0$ being an arbitrary point of $\RR^N$. But then, if $\Om$
is bounded, one gets
\begin{equation}
  \| \tmmathbf{u} \|_{L^p ( \Om )} \leqslant C_{p, N}  \left(
  \inf_{x_0 \in \RR^N} \sup_{x \in \Om} | x - x_0 | \right)  \| \gradsym
  \tmmathbf{u} \|_{L^p ( \Om )}
\end{equation}
from which \eqref{eq:PKCc3} follows as a special case. \jqed

\smallskip

{\noindent}\tmstrong{Proof of
Theorem~\ref{thm:mainPK}.{\tmem{ii.}}} {\tmname{Proof of
\eqref{eq:PKCc4}.}} If $\tmmathbf{u} \in C^{\infty} ( \overline{\Om},
\RR^N )$ and $\Om$ is a bounded Lipschitz domain, then integrating on
$\Om$ both sides of \eqref{eq:fundrel}, passing to the limit for $\varepsilon
\to 0$, and taking into account \eqref{eq:temptoest2}, we get that
\begin{equation}
  \int_{\Om} | \tmmathbf{u} |^p \leqslant C_{p, N} \int_{x \in \Om} |
  \tmmathbf{u} |^{p - 1} | \gradsym \tmmathbf{u} |  | x | + \left(
  \frac{p + 1}{p + N} \right) \int_{\xi \in \partial \Om} | \tmmathbf{u} |^p 
  | \xi | .
\end{equation}
Again, the previous computation works verbatim if one replaces everywhere $x$
with $x - x_0$, with $x_0$ an arbitrary point of $\RR^N$. Therefore, we have
\begin{equation}
  \int_{\Om} | \tmmathbf{u} |^p \leqslant ( \tmop{diam} \Om ) C_{p,
  N} \int_{\Om} | \tmmathbf{u} |^{p - 1} | \gradsym \tmmathbf{u} |
  + \left( \frac{p + 1}{p + N} \right) ( \tmop{diam} \Om )
  \int_{\partial \Om} | \tmmathbf{u} |^p  \, .\label{eq:temptoest3}
\end{equation}
Now, if $p = 1$, we are done. Otherwise, by Young's inequality for products
applied to the first term on the right-hand side of \eqref{eq:temptoest3}, we
obtain that
\begin{equation}
  \frac{1}{p} \int_{\Om} | \tmmathbf{u} |^p \leqslant \frac{( \tmop{diam}
  \Om )^p C_{p, N}^p}{p} \int_{\Om} | \gradsym \tmmathbf{u}
  |^p + \left( \frac{p + 1}{p + N} \right) ( \tmop{diam} \Om
  ) \int_{\partial \Om} | \tmmathbf{u} |^p \, ,
\end{equation}
from which \eqref{eq:PKCc4} follows.\jqed

\section*{Acknowledgments}
\noindent
{\tmname{G.Di~F.}} and {\tmname{F.S.}} are members of Gruppo Nazionale per l'Analisi Matematica, la Probabilità e le loro Applicazioni (GNAMPA) of INdAM. {\tmname{G.Di~F.}} and {\tmname{F.S.}} acknowledge support from  the Italian Ministry of Education and Research through the PRIN2022 project {\emph{Variational Analysis of Complex Systems in Material Science, Physics and Biology}} No.~2022HKBF5C.

{\tmname{G.Di~F.}} acknowledges support from the Austrian Science Fund (FWF) through the project {\emph{Analysis and Modeling of Magnetic Skyrmions}} (grant 10.55776/P34609). {\tmname{G.Di~F.}} also thanks the Hausdorff Research Institute for Mathematics in Bonn for its hospitality during the Trimester Program {\tmem{Mathematics for
Complex Materials}} funded by the Deutsche Forschungsgemeinschaft (DFG, German
Research Foundation) under Germany Excellence Strategy -- EXC-2047/1 --
390685813. {\tmname{G.Di~F.}} also thanks TU Wien and MedUni Wien for their hospitality.

The work of {\tmname{F.S.}} was partially supported by the project \emph{Variational methods for stationary and evolution problems with singularities and interfaces} PRIN 2017 (2017BTM7SN) financed by the Italian Ministry of Education, University, and Research and by the project Starplus 2020 Unina Linea~1 \emph{New challenges in the variational modeling of continuum mechanics} from the University of Naples ``Federico II'' and Compagnia di San Paolo (CUP: E65F20001630003).

\bibliographystyle{acm}
\bibliography{literature3}

\begin{thebibliography}{10}

\bibitem{Acosta2017a}
{\sc Acosta, G., and Dur\'{a}n, R.~G.}
\newblock {\em Divergence operator and related inequalities}.
\newblock SpringerBriefs in Mathematics. Springer, New York, 2017.

\bibitem{Beltrami_1968}
{\sc Beltrami, E.~J.}
\newblock Another proof of {W}eyl's lemma.
\newblock {\em {SIAM} Review 10}, 2 (apr 1968), 212--213.

\bibitem{Breit_2021}
{\sc Breit, D., and Cianchi, A.}
\newblock Symmetric gradient sobolev spaces endowed with
  rearrangement-invariant norms.
\newblock {\em Advances in Mathematics 391\/} (Nov. 2021), 107954.

\bibitem{Breit_2017}
{\sc Breit, D., Cianchi, A., and Diening, L.}
\newblock Trace-free korn inequalities in orlicz spaces.
\newblock {\em SIAM Journal on Mathematical Analysis 49}, 4 (Jan. 2017),
  2496--2526.

\bibitem{Breit2012}
{\sc Breit, D., and Diening, L.}
\newblock Sharp conditions for {K}orn inequalities in {O}rlicz spaces.
\newblock {\em Journal of Mathematical Fluid Mechanics 14}, 3 (2012), 565--573.

\bibitem{Cagnetti2022}
{\sc Cagnetti, F., Chambolle, A., and Scardia, L.}
\newblock Korn and {P}oincar\'{e}-{K}orn inequalities for functions with a
  small jump set.
\newblock {\em Mathematische Annalen 383}, 3-4 (2022), 1179--1216.

\bibitem{Chandler_Wilde_2017}
{\sc Chandler-Wilde, S.~N., Hewett, D.~P., and Moiola, A.}
\newblock Sobolev spaces on non-{L}ipschitz subsets of $\mathbb{R}^n$ with
  application to boundary integral equations on fractal screens.
\newblock {\em Integral Equations and Operator Theory 87}, 2 (feb 2017),
  179--224.

\bibitem{MR3250367}
{\sc Cianchi, A.}
\newblock Korn type inequalities in {O}rlicz spaces.
\newblock {\em Journal of Functional Analysis 267}, 7 (2014), 2313--2352.

\bibitem{Ciarlet2010}
{\sc Ciarlet, P.~G.}
\newblock On {K}orn's inequality.
\newblock {\em Chinese Annals of Mathematics. Series B 31}, 5 (2010), 607--618.

\bibitem{Ciarlet[2022]copyright2022}
{\sc Ciarlet, P.~G.}
\newblock {\em Mathematical elasticity. {V}olume {I}. {T}hree-dimensional
  elasticity}, vol.~84 of {\em Classics in Applied Mathematics}.
\newblock Society for Industrial and Applied Mathematics (SIAM), Philadelphia,
  PA, 2022.
\newblock Reprint of the 1988 edition [ 0936420].

\bibitem{Ciarlet[2022]copyright2022b}
{\sc Ciarlet, P.~G.}
\newblock {\em Mathematical elasticity. {V}olume {II}. {T}heory of plates},
  vol.~85 of {\em Classics in Applied Mathematics}.
\newblock Society for Industrial and Applied Mathematics (SIAM), Philadelphia,
  PA, 2022.
\newblock Reprint of the 1997 edition [ 1477663].

\bibitem{Ciarlet[2022]copyright2022a}
{\sc Ciarlet, P.~G.}
\newblock {\em Mathematical elasticity. {V}olume {III}. {T}heory of shells},
  vol.~86 of {\em Classics in Applied Mathematics}.
\newblock Society for Industrial and Applied Mathematics (SIAM), Philadelphia,
  PA, 2022.
\newblock Reprint of [ 1757535].

\bibitem{MR3096302}
{\sc Ciarlet, P.~G., and Mardare, S.}
\newblock Nonlinear {S}aint-{V}enant compatibility conditions and the intrinsic
  approach for nonlinearly elastic plates.
\newblock {\em Mathematical Models and Methods in Applied Sciences 23}, 12
  (2013), 2293--2321.

\bibitem{Contilecture}
{\sc Conti, S.}
\newblock {Geometric rigidity for compatible and almost-compatible fields},
  2023.
\newblock Lectures given at the "Vito Volterra Meeting", Department of
  Mathematics, Sapienza, Rome, June 2023.

\bibitem{MR3194689}
{\sc Conti, S., Dolzmann, G., and M\"{u}ller, S.}
\newblock Korn's second inequality and geometric rigidity with mixed growth
  conditions.
\newblock {\em Calculus of Variations and Partial Differential Equations 50},
  1-2 (2014), 437--454.

\bibitem{Conti2005}
{\sc Conti, S., Faraco, D., and Maggi, F.}
\newblock A new approach to counterexamples to {$L^1$} estimates: {K}orn's
  inequality, geometric rigidity, and regularity for gradients of separately
  convex functions.
\newblock {\em Archive for Rational Mechanics and Analysis 175}, 2 (2005),
  287--300.

\bibitem{Conti_2021}
{\sc Conti, S., and Garroni, A.}
\newblock Sharp rigidity estimates for incompatible fields as a consequence of
  the {B}ourgain-{B}rezis div-curl result.
\newblock {\em Comptes Rendus. Mathématique 359}, 2 (Mar. 2021), 155--160.

\bibitem{Conti_2022}
{\sc Conti, S., and Gmeineder, F.}
\newblock {$\mathscr{A}$}-quasiconvexity and partial regularity.
\newblock {\em Calculus of Variations and Partial Differential Equations 61}, 6
  (Oct. 2022).

\bibitem{Dain_2006}
{\sc Dain, S.}
\newblock Generalized {K}orn’s inequality and conformal {K}illing vectors.
\newblock {\em Calculus of Variations and Partial Differential Equations 25}, 4
  (Feb. 2006), 535--540.

\bibitem{Dautray1990}
{\sc Dautray, R., and Lions, J.-L.}
\newblock {\em Mathematical analysis and numerical methods for science and
  technology. {V}ol. 4}.
\newblock Springer-Verlag, Berlin, 1990.

\bibitem{Davoli2022a}
{\sc Davoli, E., Molchanova, A., and Stefanelli, U.}
\newblock Equilibria of charged hyperelastic solids.
\newblock {\em SIAM Journal on Mathematical Analysis 54}, 2 (2022), 1470--1487.

\bibitem{Leeuw1962}
{\sc de~Leeuw, K., and Mirkil, H.}
\newblock Majorations dans {L{$_{\infty }$}} des op\'{e}rateurs
  diff\'{e}rentiels \`a coefficients constants.
\newblock {\em Comptes Rendus Hebdomadaires des S\'{e}ances de l'Acad\'{e}mie
  des Sciences 254\/} (1962), 2286--2288.

\bibitem{Desvillettes2002}
{\sc Desvillettes, L., and Villani, C.}
\newblock On a variant of {K}orn's inequality arising in statistical mechanics.
\newblock {\em ESAIM. Control, Optimisation and Calculus of Variations.
  European Series in Applied and Industrial Mathematics 8\/} (2002), 603--619.
\newblock A tribute to J. L. Lions.

\bibitem{Di_Fratta_2020bb}
{\sc {Di~Fratta}, G., and Fiorenza, A.}
\newblock A short proof of local regularity of distributional solutions of
  {P}oisson's equation.
\newblock {\em Proceedings of the American Mathematical Society 148}, 5 (2020),
  2143--2148.

\bibitem{DiFratta2022}
{\sc {Di~Fratta}, G., and Fiorenza, A.}
\newblock A unified divergent approach to {H}ardy{\textendash}{P}oincar{\'{e}}
  inequalities in classical and variable {S}obolev spaces.
\newblock {\em Journal of Functional Analysis 283}, 5 (2022), 109552.

\bibitem{Diening_2010}
{\sc Diening, L., Ruzicka, M., and Schumacher, K.}
\newblock A decomposition technique for {J}ohn domains.
\newblock {\em Annales Academiae Scientiarum Fennicae Mathematica 35\/} (Mar.
  2010), 87--114.

\bibitem{Duvaut1972}
{\sc Duvaut, G., and Lions, J.-L.}
\newblock {\em Les in\'{e}quations en m\'{e}canique et en physique}.
\newblock Travaux et Recherches Math\'{e}matiques, No. 21. Dunod, Paris, 1972.

\bibitem{Evans2010}
{\sc Evans, L.~C.}
\newblock {\em Partial differential equations}, second~ed., vol.~19 of {\em
  Graduate Studies in Mathematics}.
\newblock American Mathematical Society, Providence, RI, 2010.

\bibitem{MR3827187}
{\sc Friedrich, M.}
\newblock A piecewise {K}orn inequality in {$SBD$} and applications to
  embedding and density results.
\newblock {\em SIAM Journal on Mathematical Analysis 50}, 4 (2018), 3842--3918.

\bibitem{Friedrichs1947}
{\sc Friedrichs, K.~O.}
\newblock {On the Boundary-Value Problems of the Theory of Elasticity and
  Korn's Inequality}.
\newblock {\em The Annals of Mathematics 48}, 2 (1947), 441--471.

\bibitem{FJM}
{\sc Friesecke, G., James, R.~D., and Müller, S.}
\newblock A theorem on geometric rigidity and the derivation of nonlinear plate
  theory from three-dimensional elasticity.
\newblock {\em Communications on Pure and Applied Mathematics 55}, 11 (aug
  2002), 1461--1506.

\bibitem{Geymonat1998}
{\sc Geymonat, G., and Gilardi, G.}
\newblock Contre-exemples \`a l'in\'{e}galit\'{e} de {K}orn et au lemme de
  {L}ions dans des domaines irr\'{e}guliers.
\newblock In {\em \'{E}quations aux d\'{e}riv\'{e}es partielles et
  applications}. Gauthier-Villars, \'{E}d. Sci. M\'{e}d. Elsevier, Paris, 1998,
  pp.~541--548.

\bibitem{Gmeineder_2023}
{\sc Gmeineder, F., Lewintan, P., and Neff, P.}
\newblock Optimal incompatible {K}orn–{M}axwell–{S}obolev inequalities in
  all dimensions.
\newblock {\em Calculus of Variations and Partial Differential Equations 62}, 6
  (June 2023).

\bibitem{Gmeineder_2024}
{\sc Gmeineder, F., Lewintan, P., and Neff, P.}
\newblock {K}orn–{M}axwell–{S}obolev inequalities for general
  incompatibilities.
\newblock {\em Mathematical Models and Methods in Applied Sciences 34}, 03
  (2024), 523--570.

\bibitem{Gmeineder_2019}
{\sc Gmeineder, F., and Raiţă, B.}
\newblock Embeddings for $\mathbb{A}$-weakly differentiable functions on
  domains.
\newblock {\em Journal of Functional Analysis 277}, 12 (Dec. 2019), 108278.

\bibitem{Gmeineder_2021}
{\sc Gmeineder, F., and Spector, D.}
\newblock On {K}orn-{M}axwell-{S}obolev inequalities.
\newblock {\em Journal of Mathematical Analysis and Applications 502}, 1 (Oct.
  2021), 125226.

\bibitem{Gobert1962}
{\sc Gobert, J.}
\newblock Une in\'{e}galit\'{e} fondamentale de la th\'{e}orie de
  l'\'{e}lasticit\'{e}.
\newblock {\em Bulletin de la Soci\'{e}t\'{e} Royale des Sciences de Li\`ege
  31\/} (1962), 182--191.

\bibitem{Grabovsky2018}
{\sc Grabovsky, Y., and Harutyunyan, D.}
\newblock Korn inequalities for shells with zero {G}aussian curvature.
\newblock {\em Annales de l'Institut Henri Poincar\'{e} C. Analyse Non
  Lin\'{e}aire 35}, 1 (2018), 267--282.

\bibitem{Horgan1995}
{\sc Horgan, C.~O.}
\newblock Korn's inequalities and their applications in continuum mechanics.
\newblock {\em SIAM Review. A Publication of the Society for Industrial and
  Applied Mathematics 37}, 4 (1995), 491--511.

\bibitem{Jiang2017}
{\sc Jiang, R., and Kauranen, A.}
\newblock Korn's inequality and {J}ohn domains.
\newblock {\em Calculus of Variations and Partial Differential Equations 56}, 4
  (July 2017), Paper No. 109, 18.

\bibitem{Kesavan_2005}
{\sc Kesavan, S.}
\newblock On {P}oincar{\'{e}}'s and {J.L.}~{L}ions' lemmas.
\newblock {\em Comptes Rendus Mathematique 340}, 1 (2005), 27--30.

\bibitem{Kondrat_ev_1988}
{\sc Kondrat{\textquotesingle}ev, V.~A., and Oleinik, O.~A.}
\newblock Boundary-value problems for the system of elasticity theory in
  unbounded domains. {K}orn's inequalities.
\newblock {\em Russian Mathematical Surveys 43}, 5 (1988), 65--119.

\bibitem{Kondratiev1989}
{\sc Kondratiev, V.~A., and Oleinik, O.~A.}
\newblock On {K}orn's inequalities.
\newblock {\em Comptes Rendus des S\'{e}ances de l'Acad\'{e}mie des Sciences.
  S\'{e}rie I. Math\'{e}matique 308}, 16 (1989), 483--487.

\bibitem{Korn1907}
{\sc Korn, A.}
\newblock Abhandlungen zur {E}lastizitätstheorie. {D}ie {E}igenschwingungen
  eines elastischen {K}örpers mit ruhender {O}berfläche.

\bibitem{korn1908solution}
{\sc Korn, A.}
\newblock Solution g\'{e}n\'{e}rale du probl\`{e}me d'\'{e}quilibre dans la
  th\'{e}orie de l'\'{e}lasticit\'{e} dans le cas o\`{u} les efforts sont
  donn\'{e}s \`{a} la surface.
\newblock {\em Annales de la Facult\'{e} des sciences de Toulouse 10\/} (1908),
  165--269.

\bibitem{korn1909einige}
{\sc Korn, A.}
\newblock {\"U}ber einige {U}ngleichungen, welche in der {T}heorie der
  elastischen und elektrischen {S}chwingungen eine {R}olle spielen.
\newblock {\em Bulletin international de l'Acad\'{e}mie de sciences de Cracovie
  3\/} (1909), 705--724.

\bibitem{KornStory1950}
{\sc Korn, T.~M.}
\newblock {\em Trailblazer to {T}elevision: {T}he Story of {A}rthur {K}orn}.
\newblock Charles Scribner's Sons, 1950.

\bibitem{Leoni_2010}
{\sc Leoni, G., Garroni, A., and Ponsiglione, M.}
\newblock Gradient theory for plasticity via homogenization of discrete
  dislocations.
\newblock {\em Journal of the European Mathematical Society 12}, 5 (Aug. 2010),
  1231--1266.

\bibitem{Magenes1958}
{\sc Magenes, E., and Stampacchia, G.}
\newblock I problemi al contorno per le equazioni differenziali di tipo
  ellittico.
\newblock {\em Annali della Scuola Normale Superiore di Pisa - Classe di
  Scienze 12}, 3 (1958), 247--358.

\bibitem{Marsden1994}
{\sc Marsden, J.~E., and Hughes, T. J.~R.}
\newblock {\em Mathematical foundations of elasticity}.
\newblock Dover Publications, Inc., New York, 1994.
\newblock Corrected reprint of the 1983 original.

\bibitem{Mikhlin1965}
{\sc Mikhlin, S.~G.}
\newblock {\em The problem of the minimum of a quadratic functional}.
\newblock Holden-Day Series in Mathematical Physics. Holden-Day, Inc., San
  Francisco, Calif.-London-Amsterdam, 1965.
\newblock Translated by A. Feinstein.

\bibitem{Mosolov1971}
{\sc Mosolov, P.~P., and Mjasnikov, V.~P.}
\newblock A proof of {K}orn's inequality.
\newblock {\em Doklady Akademii Nauk SSSR 201\/} (1971), 36--39.

\bibitem{MR3294348}
{\sc Neff, P., Pauly, D., and Witsch, K.-J.}
\newblock Poincar\'{e} meets {K}orn via {M}axwell: extending {K}orn's first
  inequality to incompatible tensor fields.
\newblock {\em Journal of Differential Equations 258}, 4 (Feb. 2015),
  1267--1302.

\bibitem{Nitsche1981}
{\sc Nitsche, J.~A.}
\newblock On {K}orn's second inequality.
\newblock {\em RAIRO Analyse Num\'{e}rique 15}, 3 (1981), 237--248.

\bibitem{Oleinik1992}
{\sc Oleinik, O.~A., Shamaev, A.~S., and Yosifian, G.~A.}
\newblock {\em Mathematical Problems in Elasticity and Homogenization}, vol.~26
  of {\em Studies in Mathematics and its Applications}.
\newblock Elsevier Science \& Technology Books, 1992.

\bibitem{Ornstein_1962}
{\sc Ornstein, D.}
\newblock A non-inequality for differential operators in the {$L_1$} norm.
\newblock {\em Archive for Rational Mechanics and Analysis 11}, 1 (1962),
  40--49.

\bibitem{Payne1961}
{\sc Payne, L.~E., and Weinberger, H.~F.}
\newblock On {K}orn's inequality.
\newblock {\em Archive for Rational Mechanics and Analysis 8\/} (1961), 89--98.

\bibitem{Reshetnyak_1970}
{\sc Reshetnyak, Y.~G.}
\newblock Estimates for certain differential operators with finite-dimensional
  kernel.
\newblock {\em Siberian Mathematical Journal 11}, 2 (1970), 315--326.

\bibitem{Simader1996}
{\sc Simader, C.~G., and Sohr, H.}
\newblock {\em The {D}irichlet problem for the {L}aplacian in bounded and
  unbounded domains}, vol.~360 of {\em Pitman Research Notes in Mathematics
  Series}.
\newblock Longman, Harlow, 1996.

\bibitem{Smith1961}
{\sc Smith, K.~T.}
\newblock Inequalities for formally positive integro-differential forms.
\newblock {\em Bulletin of the American Mathematical Society 67\/} (1961),
  368--370.

\bibitem{Temam2019}
{\sc Temam, R.}
\newblock {\em Mathematical Problems in Plasticity}.
\newblock Dover Books on Physics. Dover Publications, Incorporated, 2019.

\bibitem{NYT2016}
{\sc {The~New~York~Times}}.
\newblock {\em Photographs by {T}elegraph: {T}elevision {N}ext?}
\newblock 1907-11-24.
\newblock News article retrievable at
  www.nytimes.com/2016/03/01/science/television-history.html.

\bibitem{NYT45}
{\sc {The~New~York~Times}}.
\newblock {\em {DR.} {KORN}, {PIONEER} {IN} {RADIOPHOTO}, {DIES}; {F}ounder of
  {P}resent {S}ystems of {P}hototelegraphy and {F}acsimile {W}as 75 {R}ealized
  {P}otentialities {F}irst {T}ransatlantic {P}hoto {W}on {F}rench, {I}talian
  {P}rizes}.
\newblock 1945-12-23.
\newblock News article retrievable at
  www.nytimes.com/1945/12/23/archives/dr-korn-pioneer-in-radiophoto-dies-founder-of-present-systems-of.html.

\bibitem{Weck1994}
{\sc Weck, N.}
\newblock Local compactness for linear elasticity in irregular domains.
\newblock {\em Mathematical Methods in the Applied Sciences 17}, 2 (1994),
  107--113.

\end{thebibliography}

\end{document}